\documentclass[11pt]{amsart}

\usepackage{amsmath, amssymb, amsthm, amsrefs}
\usepackage{latexsym}
\usepackage{indentfirst}
\usepackage{graphicx}
\usepackage{placeins}
\usepackage{booktabs}
\usepackage{algorithm}
\usepackage{algorithmic}
\usepackage{multirow}

\theoremstyle{plain}

\theoremstyle{definition}

\theoremstyle{remark}

\DeclareMathOperator{\diag}{diag}

\newcommand{\eps}{\epsilon}

\newcommand{\bbm}{\begin{bmatrix}}
\newcommand{\ebm}{\end{bmatrix}}

\newcommand{\T}{\mathrm{T}}

\begin{document}

\title[Sparsifying Preconditioner for indefinite systems]{Sparsifying
  preconditioner for pseudospectral approximations of indefinite
  systems on periodic structures}

\author{Lexing Ying} 

\address{
  Department of Mathematics and Institute for Computational and Mathematical Engineering,
  Stanford University,
  Stanford, CA 94305
}

\email{lexing@math.stanford.edu}

\thanks{This work was partially supported by the National Science
  Foundation under award DMS-0846501 and the U.S. Department of
  Energy’s Advanced Scientific Computing Research program under award
  DE-FC02-13ER26134/DE-SC0009409. The author thanks Lenya Ryzhik for
  providing computing resources.}



\begin{abstract}
  This paper introduces the sparsifying preconditioner for the
  pseudospectral approximation of highly indefinite systems on
  periodic structures, which include the frequency-domain response
  problems of the Helmholtz equation and the Schr\"odinger equation as
  examples. This approach transforms the dense system of the
  pseudospectral discretization approximately into an sparse system
  via an equivalent integral reformulation and a specially-designed
  sparsifying operator. The resulting sparse system is then solved
  efficiently with sparse linear algebra algorithms and serves as a
  reasonably accurate preconditioner. When combined with standard
  iterative methods, this new preconditioner results in small
  iteration counts. Numerical results are provided for the Helmholtz
  equation and the Schr\"odinger in both 2D and 3D to demonstrate the
  effectiveness of this new preconditioner.
\end{abstract}

\maketitle

%-----------------------------------
\section{Introduction}

This paper is concerned with the numerical solution of highly
indefinite systems on periodic structures. One example comes from the
study of the propagation of high frequency acoustic and
electromagnetic waves in periodic media, which can be described in its
simplest form by the Helmholtz equation with the periodic boundary
condition
\begin{equation}
  \left(-\Delta - \frac{\omega^2}{c(x)^2}\right) u(x) = f(x),\quad x\in \T^d := [0,1)^d,
    \label{eq:helm}
\end{equation}
where $\omega$ is the wave frequency and $c(x)$ is a periodic velocity
field. This system is highly indefinite for large values of
$\omega$. A second example is the Schr\"odinger equation with the
periodic boundary condition
\[
(-\Delta + V(x) - E) u(x) = f(x),\quad x\in [0,\ell)^d,
\]
where $V(x)$ is the potential field, $E$ is the energy level, and
$\ell$ is the system size. The solution of this system appears as an
essential step of the electronic structure calculation of quantum
many-particle systems. Typically we are interested in the regime of
large values of $\ell$ and it is convenient to rescale the system to
the unit cube via the transformation $x = \ell y$:
\begin{equation}
  (-\Delta + \ell^2 V(\ell y) - \ell^2 E) u(y) = \ell^2 f(y), \quad y\in \T^d := [0,1)^d,
    \label{eq:schr}
\end{equation}

Since the domain is compact, the operators \eqref{eq:helm} and
\eqref{eq:schr} can be non-invertible for certain values of $\omega$
and $E$, respectively. In this paper, we assume that the systems are
invertible and we are interested in the efficient and accurate
numerical solutions of these systems.

Numerical solution of \eqref{eq:helm} and \eqref{eq:schr} has been a
long-standing challenge for several reasons. First, these problems can
be almost non-invertible if $\omega$ or $E$ is a (generalized)
eigenvalue of the system. Second, the systems are highly indefinite,
i.e., the operators in these equations have large number of positive
and negative eigenvalues. Third, since the solutions of these
equations are always highly oscillatory, an accurate approximation of
the solution typically requires a large number of unknowns due to the
Nyquist theorem.

The simplest numerical approach for \eqref{eq:helm} and
\eqref{eq:schr} is probably the standard finite difference and finite
element methods. These methods result in sparse linear systems with
local stencil, thus enabling the use of sparse direct solvers such as
the nested dissection method \cite{george-1973} and the multifrontal
method \cite{duff-1983}. However, these methods typically gave wrong
dispersion relationships for high-frequency/high-energy problems and
thus fail to provide accurate solutions. One solution is to use higher
order finite difference stencils that provide more accurate dispersion
relationships. However, this comes at a price of increasing the
stencil support, which quickly makes it impossible to use the sparse
direct solvers.

One practical approach for \eqref{eq:helm} and \eqref{eq:schr} is the
spectral element methods \cite{patera-1984,canuto-2007}, which
typically use local higher order polynomial bases, such as Chebyshev
functions, within each rectangular element. These methods allow for
efficient solution with the sparse direct solvers. When the polynomial
degree is sufficient high, the spectral element methods can capture
the dispersion relationship accurately. On the other hand, their
implementations typically require much more effort.

Because of the periodic domain, the pseudospectral method
\cite{orszag-1969,gottlieb-1977,trefethen-2000} with Fourier (plane
wave) bases is highly popular for the problems \eqref{eq:helm} and
\eqref{eq:schr}. They are simple to implement and typically require a
minimum number of unknowns for a fixed accuracy among all methods
discussed above. Therefore, the pseudospectral method is arguably the most
widely used approach in engineering and industrial studies of the
systems. However, the main drawback of the pseudospectral method is
that the resulting discrete systems are dense and hence it is usually
impossible to apply the efficient sparse direct solvers. This is
indeed what this paper aims to address.

In this paper, we introduce the sparsifying preconditioner for the
pseudospectral approximation of \eqref{eq:helm} and \eqref{eq:schr}
for periodic domains. The main idea is to introduce an equivalent
integral formulation and transform the dense system numerically into a
sparse one following the idea from \cite{ying-2014}. The approximate
sparse system is then solved efficiently with the help of sparse
direct solvers and serves as a reasonably accurate preconditioner for
the original pseudospectral system. When combined with standard
iterative algorithms such as GMRES \cite{saad-1986}, this new
preconditioner results in small iteration counts.

The rest of the paper is organized as follows. Section 2 introduces
the sparsifying preconditioner after discussing the pseudospectral
discretization. Numerical results for both 2D and 3D problems are
provided in Section 3 and finally future work is discussed in Section 4.

%-----------------------------------
\section{The sparsifying preconditioner}

\subsection{The pseudospectral approximation}

Since our approach treats \eqref{eq:helm} and \eqref{eq:schr} in the
same way, it is convenient to introduce a general system for both
cases:
\begin{equation}
  (-\Delta - s + q(x)) u(x) = f(x),\quad x\in \T^d := [0,1)^d,
    \label{eq:Dsq}
\end{equation}
where $s$ is a constant shift and $q(x)$ is the inhomogeneous term.
For \eqref{eq:helm}, $s$ and $q(x)$ are given by
\[
s = \int_{\T^d} \frac{\omega^2}{c(x)^2} dx,\quad q(x) = -\frac{\omega^2}{c(x)^2}+s.
\]
For \eqref{eq:schr}, they are equal to
\[
s = \int_{\T^d} (-\ell^2 V(\ell y)+\ell^2 E) dx,\quad q(x) = \ell^2 V(\ell y)-\ell^2 E+s.
\]

The pseudospectral method discretizes the domain $\T^d=[0,1)^d$
  uniformly with a uniform grid of size $n$ in each dimension.  The
  step size $h=1/n$ is chosen to ensure that there are at least 3 to 4
  points for the typical oscillation of the solution. The grid points
  are indexed by a set
\[
J =\{(j_1,\ldots,j_d): 0\le j_1,\ldots,j_d < n\}.
\]
For each $j\in J$, we define
\[
f_j = f(jh),\quad q_j = q(jh)
\]
and let $u_i$ be the numerical approximation to $u(ih)$ to be
determined.  We also introduce a grid in the Fourier domain
\[
K = \{(k_1,\ldots,k_d): -n/2 \le k_1,\ldots,k_d < n/2\}.
\]
The forward and inverse Fourier operators $F$ and $F^{-1}$ are defined
by
\begin{align*}
  (F    f)_k &= \frac{1}{n^{d/2}} \sum_{j\in J} e^{-2\pi i (j\cdot k)/n} f_j, \quad k\in K\\
  (F^{-1}g)_j &= \frac{1}{n^{d/2}} \sum_{k\in K} e^{+2\pi i (j\cdot k)/n} g_k, \quad j\in J.
\end{align*}
The pseudospectral method discretizes the Laplacian operator with
\[
L := F^{-1} \diag(4\pi^2 |k|^2)_{k\in K} F.
\]
By a slight abuse of notation, we use $u$ to denote the vector with
entries $u_j$ for $j\in J$ and similarly for the vectors $f$ and
$q$. The discretized system of \eqref{eq:Dsq} then becomes
\begin{equation}
  (L-s+q) u = f,
  \label{eq:Lsq}
\end{equation}
where $q$ also stands for the diagonal operator of entry-wise
multiplication with the elements of the vector $q$.

\subsection{Main idea}
We assume without loss of generality that $L-s$ is invertible, which
can be easily satisfied by perturbing $s$ slightly if necessary. We
define
\begin{equation}
  G := (L-s)^{-1} = F^{-1} \diag\left(\frac{1}{4\pi^2 |k|^2-s}\right)_{k\in K} F,
  \label{eq:G}
\end{equation}
which is a discrete convolution operator that can be applied
efficiently with the fast Fourier transform. Applying $G$ to both
sides of \eqref{eq:Lsq} gives
\[
(I+ Gq) u = Gf =: g.
\] 
The main idea of the sparsifying preconditioner is to introduce a
sparse matrix $Q$ such that in the preconditioned system
\[
Q(I+Gq) u = Qg
\]
the operator $QG$ is approximately sparse as well. Based on this, we
define the matrix $P$ to be the truncated version of $Q(I+Gq)$ and
arrive at the approximate equation
\[
P u \approx Qg.
\]
Since $P$ is sparse, we factorize it with sparse direct solvers such
as the nested dissection algorithm and set
\[
u \Leftarrow P^{-1}Qg
\]
as the approximate inverse.

More precisely, for a given point $j\in J$ we denote by $\mu(j)$ its
neighborhood (to be defined below). The row $Q(j,:)$ should satisfy
the following two conditions:
\begin{itemize}
\item $Q(j,:)$ is supported on $\mu(j)$,
\item $Q(j,\mu(j)) G(\mu(j),\mu(j)^c) \approx 0$.
\end{itemize}
These two conditions imply that $Q(j,:)G(:,:) = Q(j,\mu(j))
G(\mu(j),:)$ is essentially supported in $\mu(j)$. We then define the
matrix $C$ such that each row $C(j,:)$ is supported only in $\mu(j)$
and
\[
C(j,\mu(j)) = Q(j,\mu(j)) G(\mu(j),\mu(j)).
\]
This definition implies that $C$ has the same sparsity pattern as $Q$,
$C\approx QG$, and also
\[
Q(I + G q) \approx Q + C q =: P.
\]
Since $q$ is diagonal, $P=Q+Cq$ have the same non-zero pattern as $Q$.

\subsection{Details}

There are two problems that remain to be addressed. The first one is
the definition of the neighborhood $\mu(j)$ of a point $j\in J$. For a
given set of grid points $s$, we define
\[
\gamma(s) = \{i | \exists j\in s, \|j-i\|_\infty \le 1\},
\]
where the distance is measured modulus the grid size $n$ in each
dimension. In \cite{ying-2014}, $\mu(j) = \gamma(\{j\})$, i.e.,
$\mu(j)$ contains the nearest neighbors of $j$ in the $\ell_\infty$
norm. For the matrix $C$, this corresponds to setting the elements in
$(QR)(j,\gamma(\{j\})^c)$ to zero. However, the error introduced turns
out to be too large in the current setting, since the system
considered here can be very ill-conditioned. Therefore, one needs to
increase $\mu(j)$ in order to take more off-diagonal entries into
consideration. However, because $\mu(j)$ controls the sparsity pattern
of the matrix $P$ that is to be factorized with sparse direct solvers,
an increase of $\mu(j)$ should be done in a way not to sacrifice the
efficiency of the sparse direct solvers.

\begin{figure}[h!]
  \includegraphics[height=2.5in]{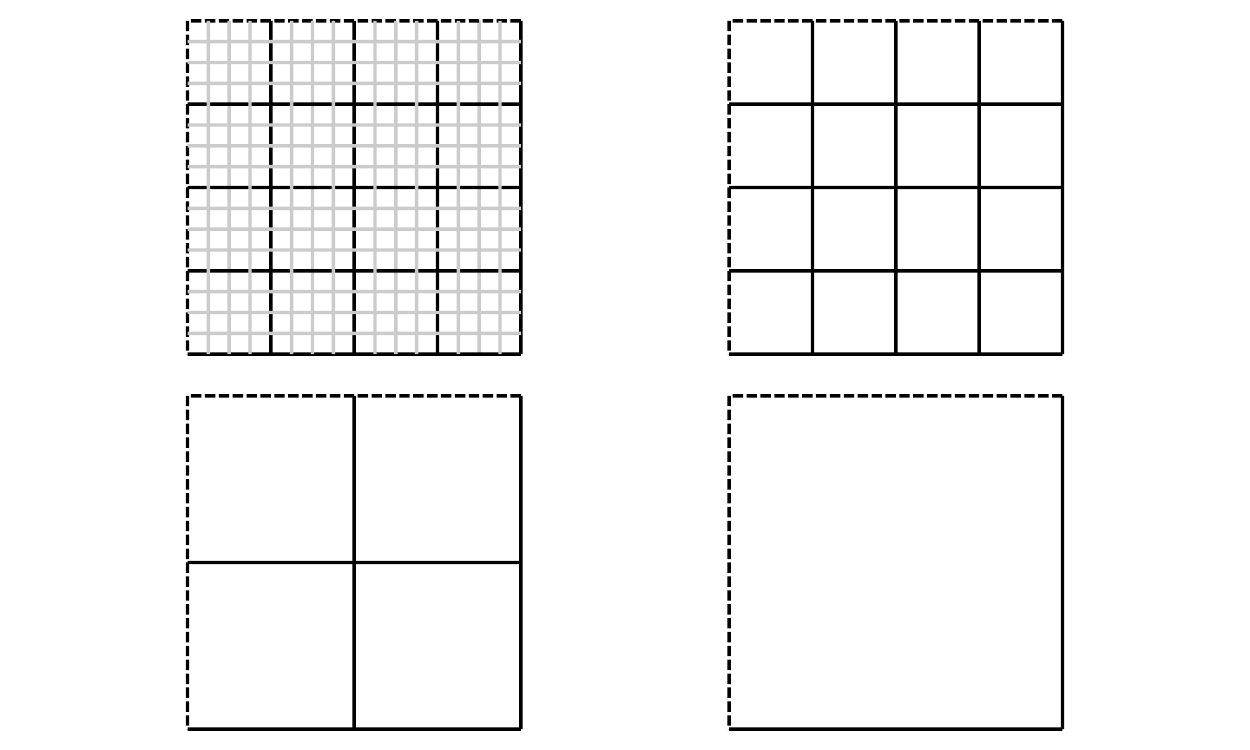}
  \caption{The nested dissection algorithm. The domain is partitioned
    recursively into smaller square boxes until the side length is
    equal to a constant $bh$. The algorithm recursively eliminates the
    interior nodes of a box via Schur complement, combines four boxes
    into a large one, and repeats the process. The dotted lines stand
    for the degrees of freedom repeated on the other side of the
    domain due to periodicity.}
  \label{fig:nd}
\end{figure}

Throughout the rest of the paper, we use the nested dissection
algorithm to factorize $P$. In this algorithm, the domain $\T^d$ is
partitioned recursively into square boxes until a certain size $bh$ is
reached in each dimension. Each leaf box contains $(b-1)^d$ points in
its interior and $(b+1)^d$ points in its closure. The algorithm then
recursively eliminates the interior nodes of a box via Schur
complement, combine four boxes into a large one, and repeat. Figure
\ref{fig:nd} illustrates the nested dissection algorithm in a 2D
setting.

An essential observation is that $\mu(j)$ can include more points
without affecting much the complexity of the nested dissection
algorithm. Let us discuss the 2D case first. At the leaf level, the
set $J$ is partitioned into a disjoint union of three types of sets:
\begin{itemize}
\item a {\em cell} set that contains the $(b-1)^2$ interior grid
  points of a leaf box,
\item an {\em edge} set that contains the $(b-1)$ grid points at the
  interface of two adjacent leaf boxes, and
\item a {\em vertex} set that contains only $1$ grid point at the
  interface of four adjacent leaf boxes.
\end{itemize}
The neighborhood $\mu(j)$ of a point $j\in J$ is defined as follows
based on the type of the set that contains it:
\begin{itemize}
\item for $j$ in a cell set $c$, we set $\mu(j) = \gamma(c)$;
\item for $j$ in an edge set $e$, we set $\mu(j) = \gamma(e)$;
\item for $j$ in an vertex set $v=\{j\}$, we set $\mu(j) = \gamma(v)$.
\end{itemize}

In the 3D case, the set $J$ is partitioned into a disjoint union of
four types of sets:
\begin{itemize}
\item a {\em cell} set that contains the $(b-1)^3$ interior grid
  points of a leaf box,
\item a {\em face} set that contains the $(b-1)^2$ grid points at the
  interface of two adjacent leaf boxes,
\item an {\em edge} set that contains the $(b-1)$ grid points at the
  interface of four adjacent leaf boxes, and
\item a {\em vertex} set that contains only $1$ grid point at the
  interface of eight adjacent leaf boxes.
\end{itemize}
For the extra case of $j$ in a face set $f$, we define $\mu(j) =
\gamma(f)$.

The second problem is the computation of $Q(j,\mu(j))$. Our choice of
the neighborhood shows that $\mu(j)$ is the same for all $j$ in the
same (cell, face, edge, or vertex) set. Therefore it is convenient to
consider all such $j$ together.

We first fix a cell set $c$ and consider all points $j\in c$. The
sparsity condition on $Q$ can be rewritten as
\begin{equation}
  Q(c,\gamma(c)) G(\gamma(c),\gamma(c)^c) \approx 0.
  \label{eq:Qc}
\end{equation}
The rows of the submatrix $Q(c,\gamma(c))$ should also be linearly
independent in order for $Q$ to be non-singular. We define $\beta(c) =
\gamma(c)\setminus c$, i.e., the set of points that are on the
boundary of $c$. Since the matrix $G$ defined in \eqref{eq:G} is the
Green's function of a discretized partial differential operator, the
rows of $G(c,\gamma(c)^c)$ can be approximated accurately using the
linear combinations of the rows of $G(\beta(c),\gamma(c)^c)$, i.e.,
there exists a matrix $T_c$ such that
\[
G(c,\gamma(c)^c) \approx T_c G(\beta(c),\gamma(c)^c).
\]
$T_c$ matrix can be computed by
\[
T_c = G(c,\gamma(c)^c) (G(\beta(c),\gamma(c)^c))^+,
\]
where $(\cdot)^+$ stands for the pseudoinverse. Finally, we define
\[
Q(c,\gamma(c)) = \begin{bmatrix} I & -T_c \end{bmatrix},
\]
assuming the columns are ordered as $(c,\beta(c))$. This submatrix
meets the condition \eqref{eq:Qc} and clearly has linearly independent
rows. We remark that the computation of $T_c$ is the same for any cell
$c$ and hence we only need to compute it once.

For a fixed face set $f$, we compute
\[
T_f = G(f,\gamma(f)^c) (G(\beta(f),\gamma(f)^c))^+
\]
with $\beta(f) = \gamma(f)\setminus f$ and set
\[
Q(f,\gamma(f)) = \begin{bmatrix} I & -T_f \end{bmatrix},
\]
assuming the columns are ordered as $(f,\beta(f))$.

For a fixed edge set $e$, we define
\[
T_e = G(e,\gamma(e)^c) (G(\beta(e),\gamma(e)^c))^+
\]
with $\beta(e) = \gamma(e)\setminus e$ and set
\[
Q(e,\gamma(e)) = \begin{bmatrix} I & -T_e \end{bmatrix},
\]
assuming the columns are ordered as $(e,\beta(e))$.

Finally, for a fixed vertex set $v$, we let
\[
T_v = G(v,\gamma(v)^c) (G(\beta(v),\gamma(v)^c))^+
\]
with $\beta(v) = \gamma(v)\setminus v$ and set
\[
Q(v,\gamma(v)) = \begin{bmatrix} I & -T_v \end{bmatrix},
\]
assuming the columns are ordered as $(v,\beta(v))$.

Once the matrix $Q$ has been constructed, the matrix $C$ is computed
as follows. For a cell set $c$, we set
\[
C(c,\gamma(c)) = Q(c,\gamma(c)) G(\gamma(c),\gamma(c)),
\]
and similarly for a face set $f$, an edge set $e$, and a vertex set
$v$. We recall that both $C$ and $P = Q + Cq$ have the same sparsity
pattern as $Q$.

\subsection{Complexity}

We now analyze the complexity of constructing and applying the
sparsifying preconditioner. Let $bh$ be the width of the leaf box of
the nested dissection algorithm.

In 2D, the construction algorithm consists of two parts: (i) computing
the pseudoinverses while forming $T_c$, $T_e$, and $T_v$, and (ii)
building the nested dissection factorization for $P$. The former takes
at most $O(b^4 n^2) = O(b^4 N)$ steps, while the latter takes $O(n^3 +
b^6(n/b)^2) = O(N^{3/2} + b^4 N)$ steps. Therefore, the overall
complexity of the construction algorithm is $O(N^{3/2} + b^4 N)$. The
application algorithm is essentially a nested dissection solve, which
costs $O(n^2 \log n + b^4 (n/b)^2) = O(N \log N + b^2 N)$ steps.

In 3D, the construction algorithm again consists of the same two
parts. The pseudoinverses cost $O(b^6 n^3) = O(b^6 N)$ steps, while
the nest dissection factorization takes $O(n^6 + b^9 (n/b)^3) = O(N^2
+ b^6 N)$ steps. Hence, the overall construction cost is $O(N^2 + b^6
N)$.  The application cost is equal to $O(n^4 + b^6(n/b)^3) =
O(N^{4/3} + b^3 N)$ due to a 3D nested dissection solve.

There is a clear trade-off for the choice of $b$. For small values of
$b$, the preconditioner is less efficient due to the small support of
$Q$, while the computational complexity is low. On the other hand, for
large values of $b$, the preconditioner is more effective, but the
cost is higher. In our numerical results, we set $b = O(n^{1/2})$. For
this choice, the construction and application costs in 2D are $O(N^2)$
and $O(N^{3/2})$, respectively. In 3D, they are $O(N^2)$ and
$O(N^{3/2})$ as well, respectively.

%-----------------------------------
\section{Numerical results}

The sparsifying preconditioner, as well as the nested dissection
algorithm involved, is implemented in Matlab. The numerical results
are obtained on a Linux computer with CPU speed at 2.0GHz. The GMRES
algorithm is used as the iterative solver with the relative tolerance
equal to $10^{-6}$.

\subsection{Helmholtz equation}
\begin{figure}[h!]
  \begin{center}
    \begin{tabular}{|ccc|cc|cc|}
      \hline
      $\omega/(2\pi)$ & $N$ & $b$ & $T_s$(sec) & $T_a$(sec) & $n_p$ & $T_p$(sec) \\
      \hline
      16 & $48^2$ & 3 & 2.5e-01 & 3.7e-02 & 7.0e+00 & 2.1e-01\\
      32 & $96^2$ & 6 & 5.5e-01 & 2.9e-02 & 8.0e+00 & 2.6e-01\\
      64 & $192^2$ & 6 & 2.5e+00 & 8.1e-02 & 1.1e+01 & 1.0e+00\\
      128 & $384^2$ & 12 & 1.2e+01 & 2.9e-01 & 2.8e+01 & 9.8e+00\\
      \hline
    \end{tabular}
    \includegraphics[height=1.8in]{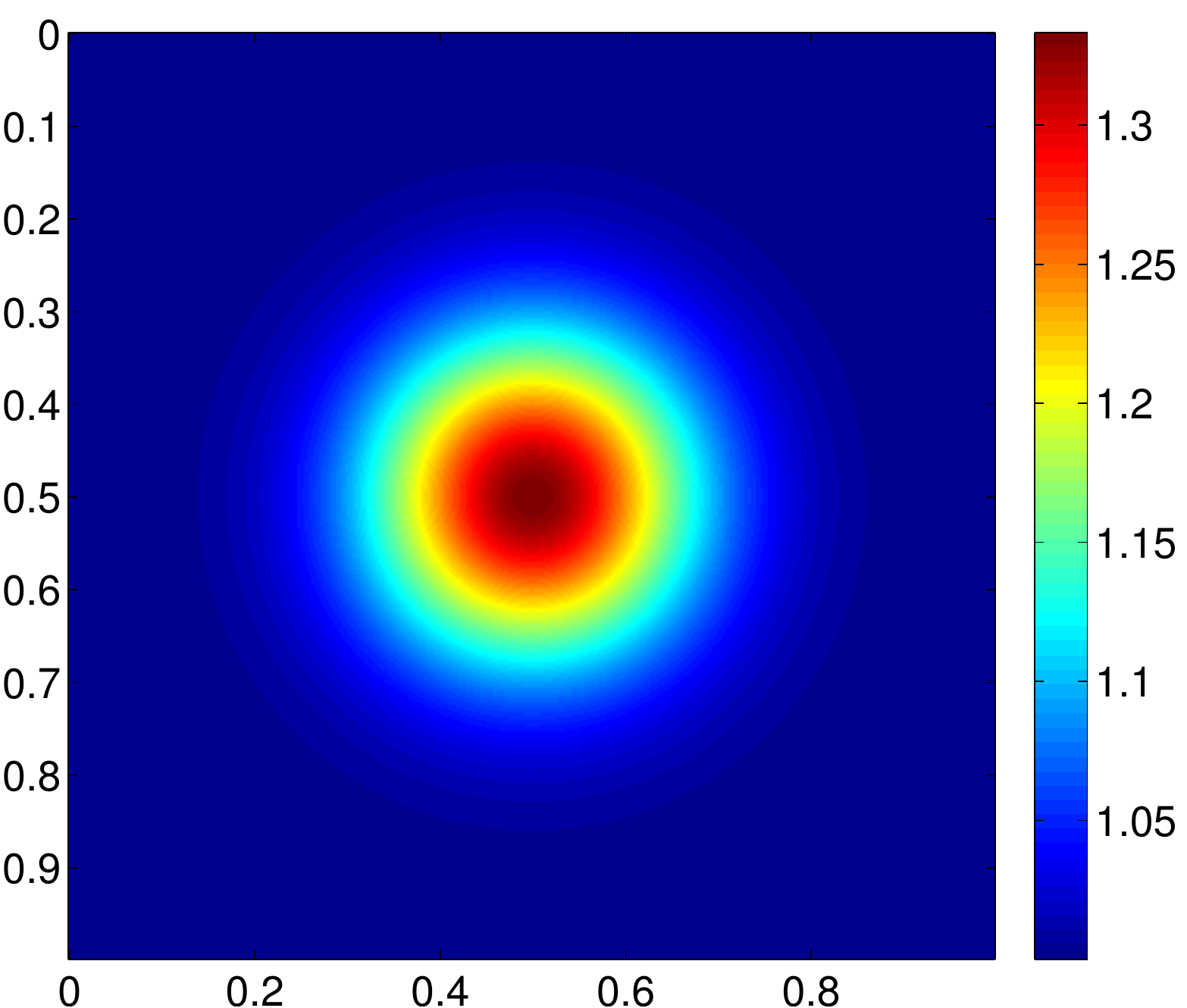} \hspace{0.25in} \includegraphics[height=1.8in]{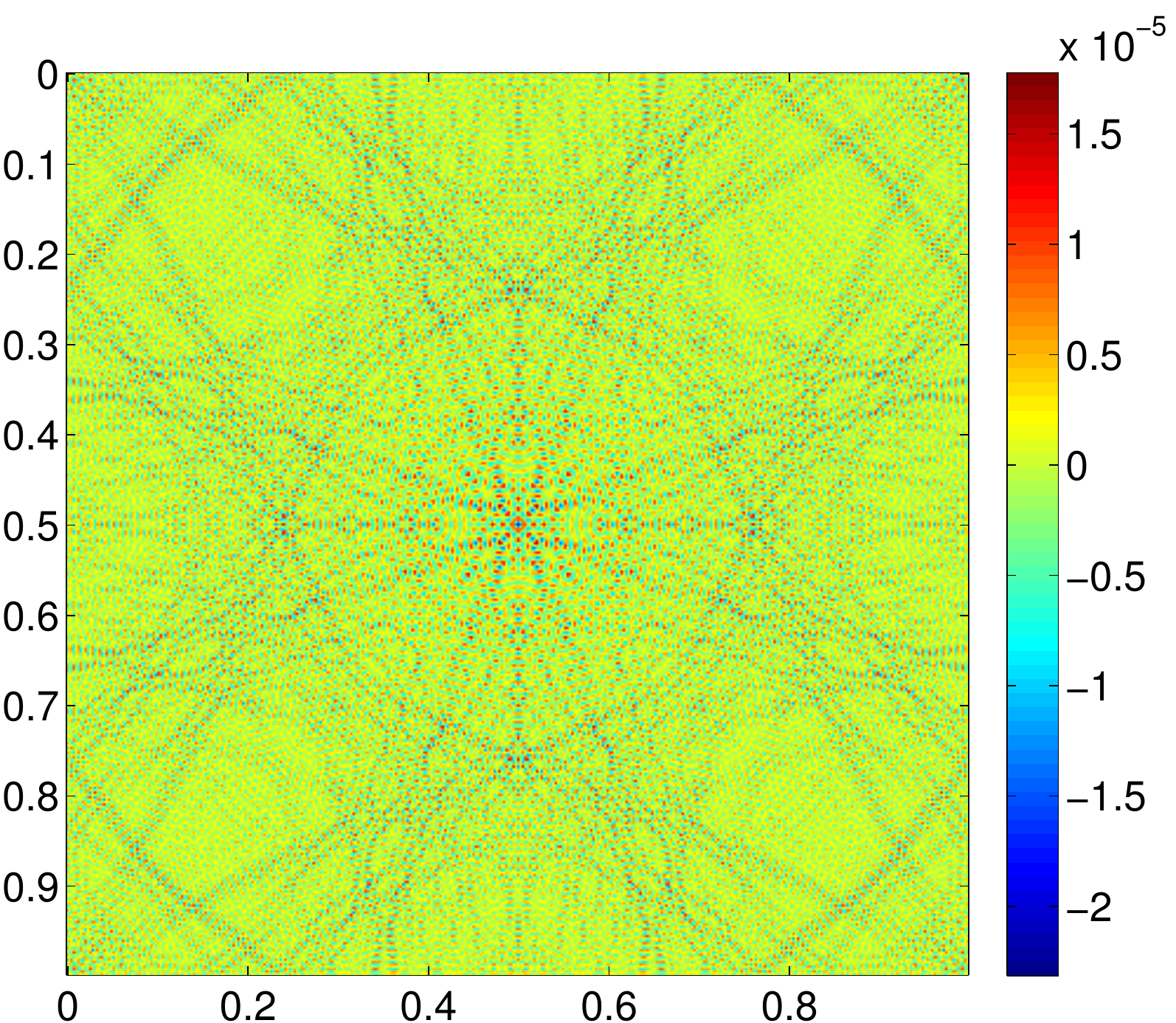}\\
  \end{center}
  \caption{Example 1 of the 2D Helmholtz equation. Top: numerical
    results.  Bottom: $c(x)$ (left) and $u(x)$ (right) for the
    largest $\omega$ value.}
  \label{fig:W21}
\end{figure}

\begin{figure}[h!]
  \begin{center}
    \begin{tabular}{|ccc|cc|cc|}
      \hline
      $\omega/(2\pi)$ & $N$ & $b$ & $T_s$(sec) & $T_a$(sec) & $n_p$ & $T_p$(sec) \\
      \hline
      16 & $48^2$ & 3 & 2.2e-01 & 1.3e-02 & 8.0e+00 & 1.1e-01\\
      32 & $96^2$ & 6 & 5.0e-01 & 2.7e-02 & 1.1e+01 & 3.8e-01\\
      64 & $192^2$ & 6 & 3.0e+00 & 7.7e-02 & 2.1e+01 & 2.1e+00\\
      128 & $384^2$ & 12 & 1.2e+01 & 2.9e-01 & 3.8e+01 & 1.3e+01\\
      \hline
    \end{tabular}
    \includegraphics[height=1.8in]{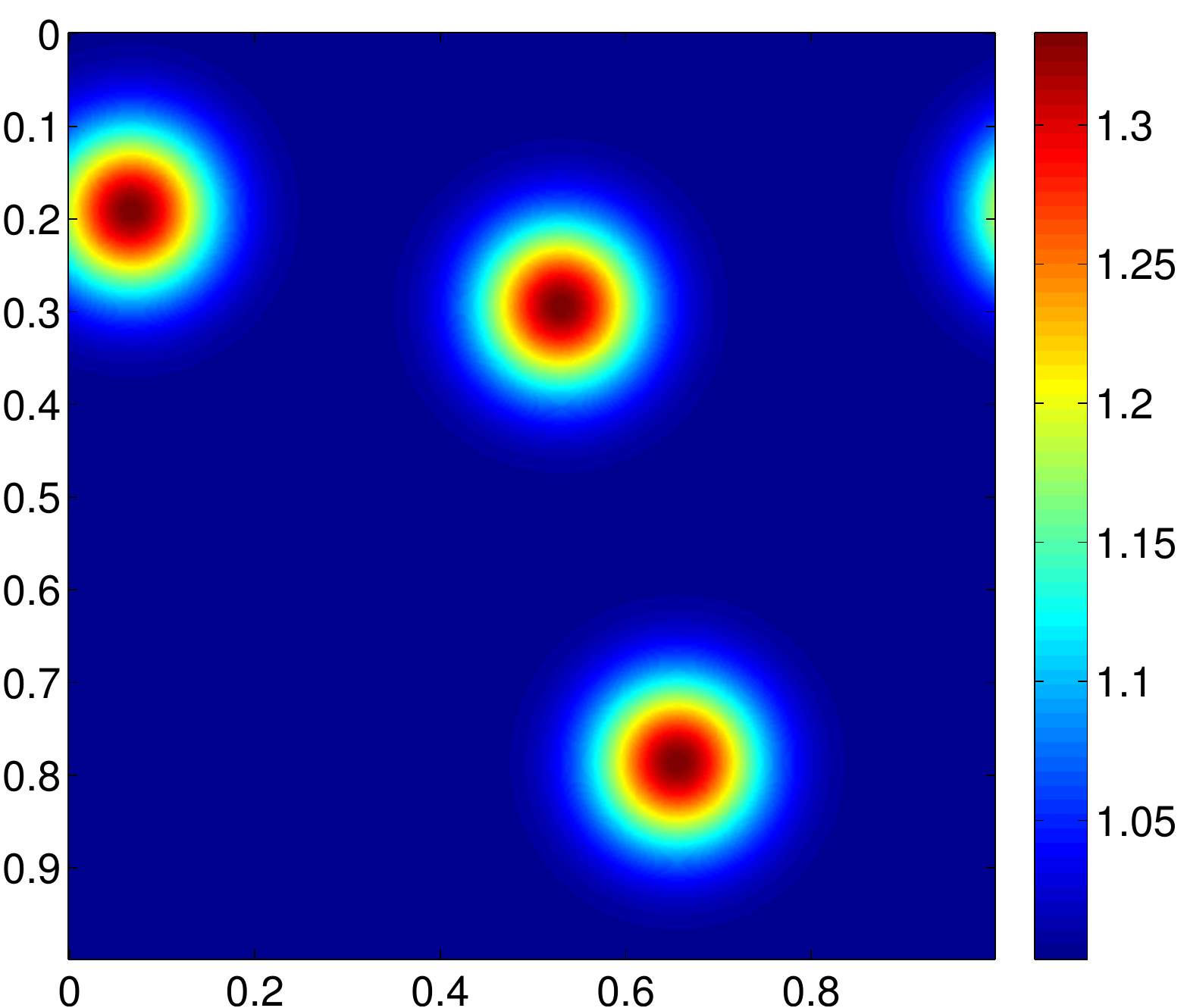} \hspace{0.25in} \includegraphics[height=1.8in]{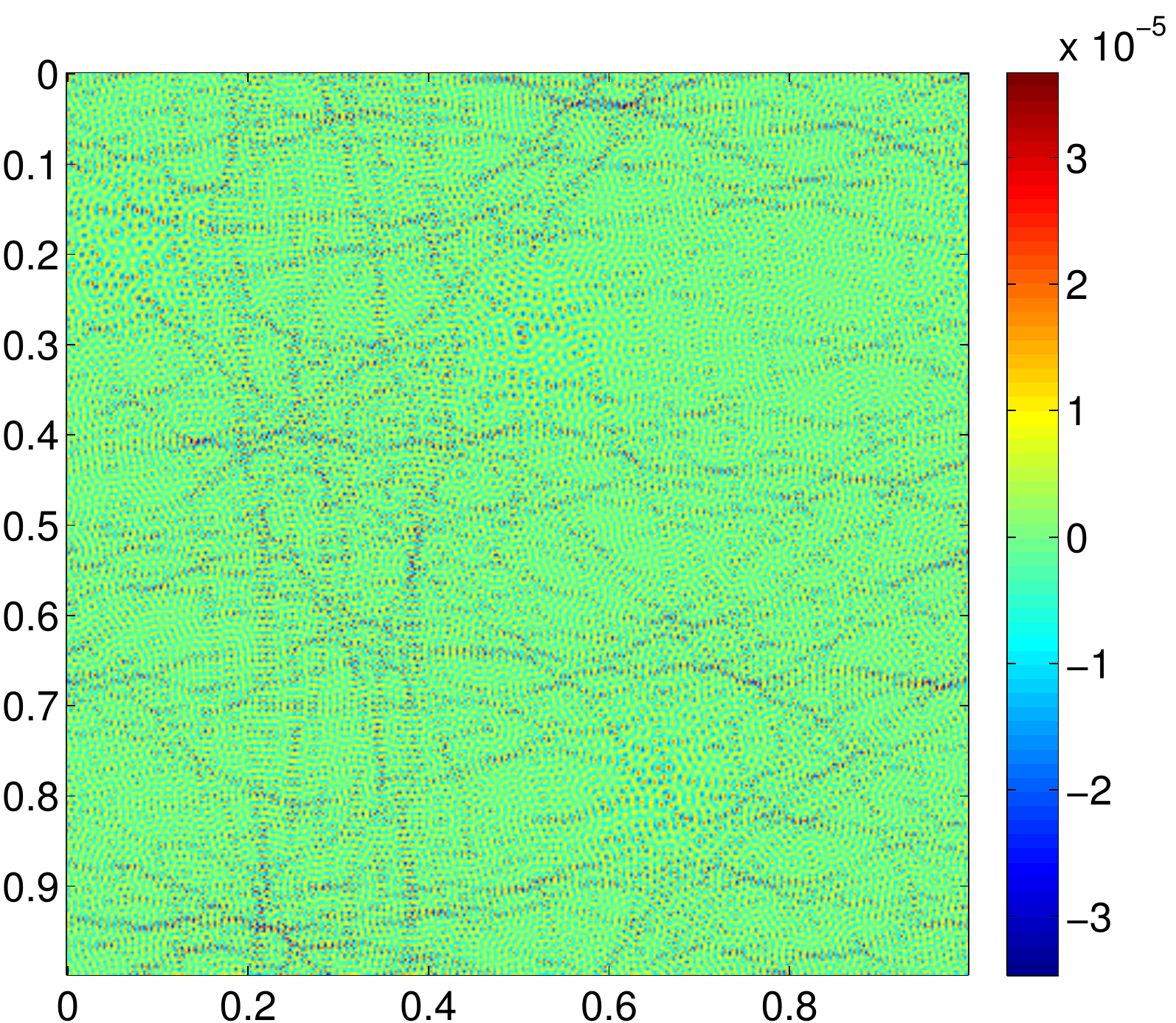}\\
  \end{center}
  \caption{Example 2 of the 2D Helmholtz equation. Top: numerical
    results.  Bottom: $c(x)$ (left) and $u(x)$ (right) for the
    largest $\omega$ value.}
  \label{fig:W22}
\end{figure}

Two tests are performed for the 2D Helmholtz equation. In the first
one the velocity field $c(x)$ is equal to one plus a Gaussian function
at the domain center, while in the second test $c(x)$ is given by one
plus three randomly placed Gaussian functions. In both tests, the
right hand side is a delta source at the center of the domain. The
results of these two tests are summarized in Figures \ref{fig:W21} and
\ref{fig:W22}. The columns of the tables are listed as follows:
\begin{itemize}
\item $\omega/(2\pi)$ is the wave number (roughly the number of oscillation across the domain),
\item $N$ is the number of unknowns,
\item $b$ is the ratio between the width of the leaf box and the step size $h$,
\item $T_s$ is the setup time of the preconditioner in seconds,
\item $T_a$ is the application time of the preconditioner in seconds,
\item $n_p$ is the iteration number of the preconditioned iteration, and
\item $T_p$ is the solution time of the preconditioned iteration in seconds.
\end{itemize}

\begin{figure}[h!]
  \begin{center}
    \begin{tabular}{|ccc|cc|cc|}
      \hline
      $\omega/(2\pi)$ & $N$ & $b$ & $T_s$(sec) & $T_a$(sec) & $n_p$ & $T_p$(sec) \\
      \hline
      4 & $12^3$ & 3 & 4.8e-01 & 1.9e-02 & 5.0e+00 & 1.7e-01\\
      8 & $24^3$ & 6 & 6.0e+00 & 4.7e-02 & 7.0e+00 & 3.9e-01\\
      16 & $48^3$ & 6 & 1.5e+02 & 4.8e-01 & 8.0e+00 & 4.0e+00\\
      32 & $96^3$ & 12 & 4.8e+03 & 7.4e+00 & 1.0e+01 & 8.4e+01\\
      \hline
    \end{tabular}
    \includegraphics[height=1.8in]{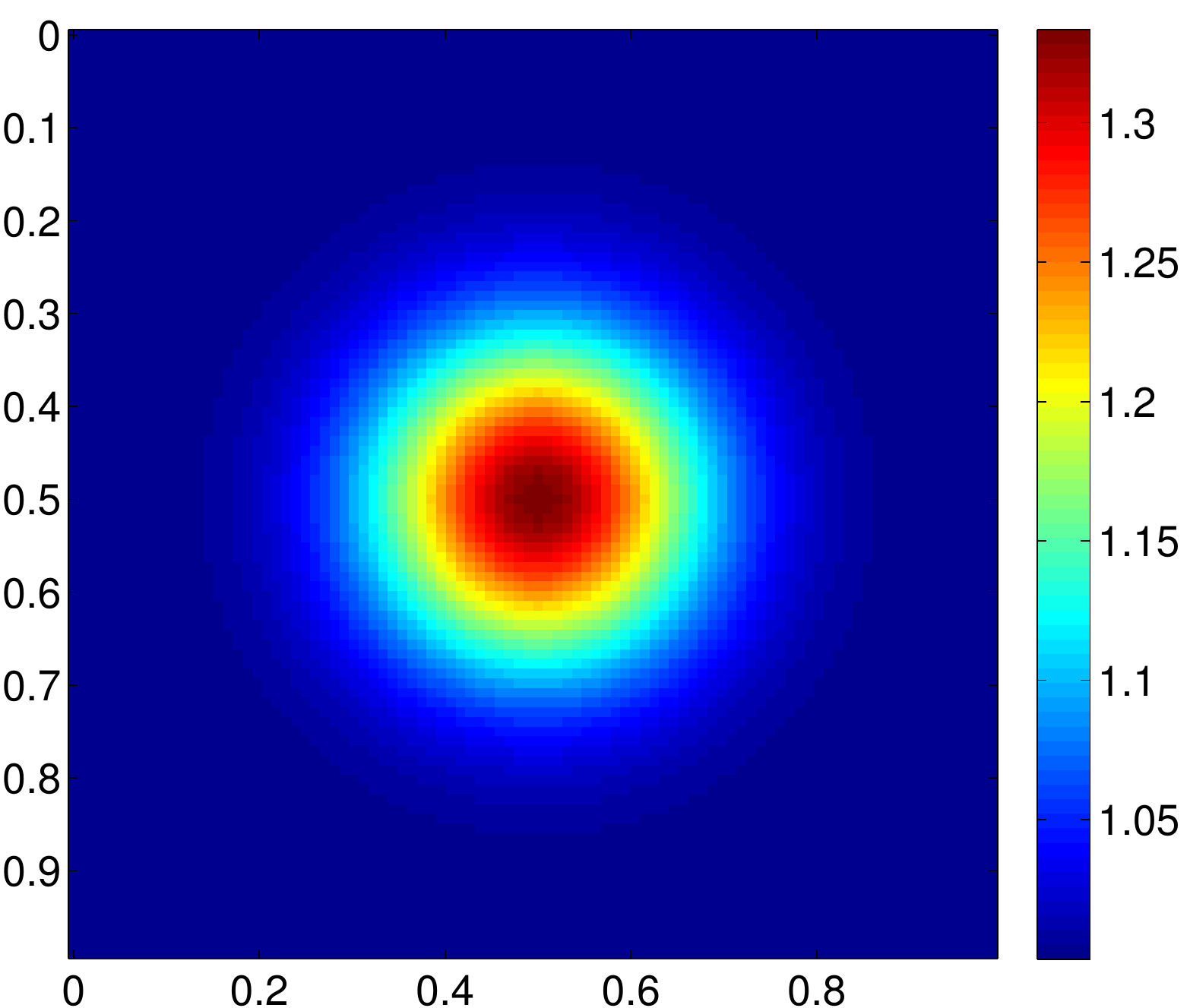} \hspace{0.25in} \includegraphics[height=1.8in]{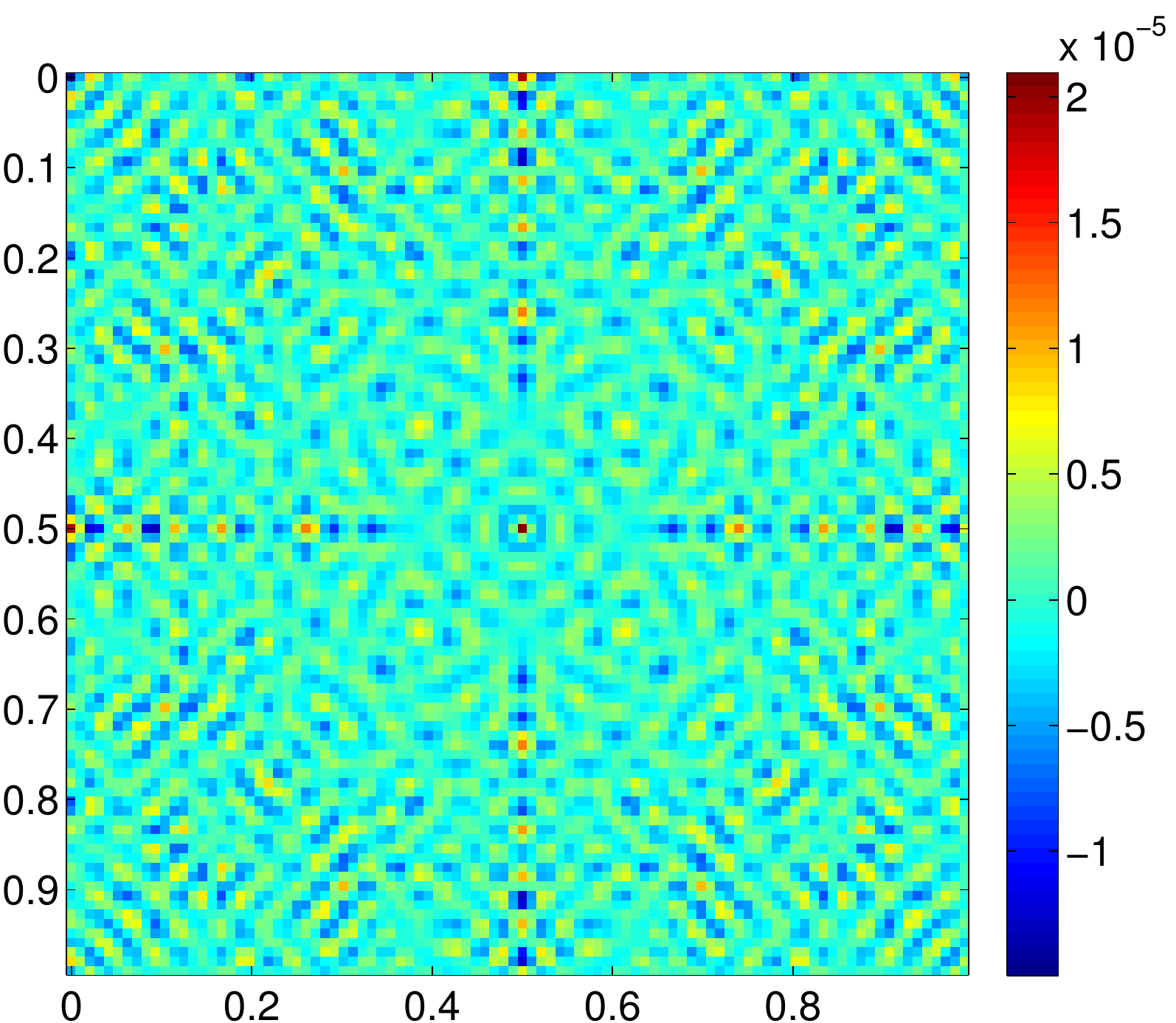}\\
  \end{center}
  \caption{Example 1 of the 3D Helmholtz equation. Top: numerical
    results.  Bottom: $c(x)$ (left) and $u(x)$ (right) for the
    largest $\omega$ value at the middle slice.}
  \label{fig:W31}
\end{figure}

\begin{figure}[h!]
  \begin{center}
    \begin{tabular}{|ccc|cc|cc|}
      \hline
      $\omega/(2\pi)$ & $N$ & $b$ & $T_s$(sec) & $T_a$(sec) & $n_p$ & $T_p$(sec) \\
      \hline
      4 & $12^3$ & 3 & 3.0e-01 & 6.5e-03 & 4.0e+00 & 2.0e-02\\
      8 & $24^3$ & 6 & 6.3e+00 & 5.1e-02 & 5.0e+00 & 2.9e-01\\
      16 & $48^3$ & 6 & 1.5e+02 & 4.8e-01 & 9.0e+00 & 4.6e+00\\
      32 & $96^3$ & 12 & 4.7e+03 & 8.4e+00 & 1.8e+01 & 1.6e+02\\
      \hline
    \end{tabular}
    \includegraphics[height=1.8in]{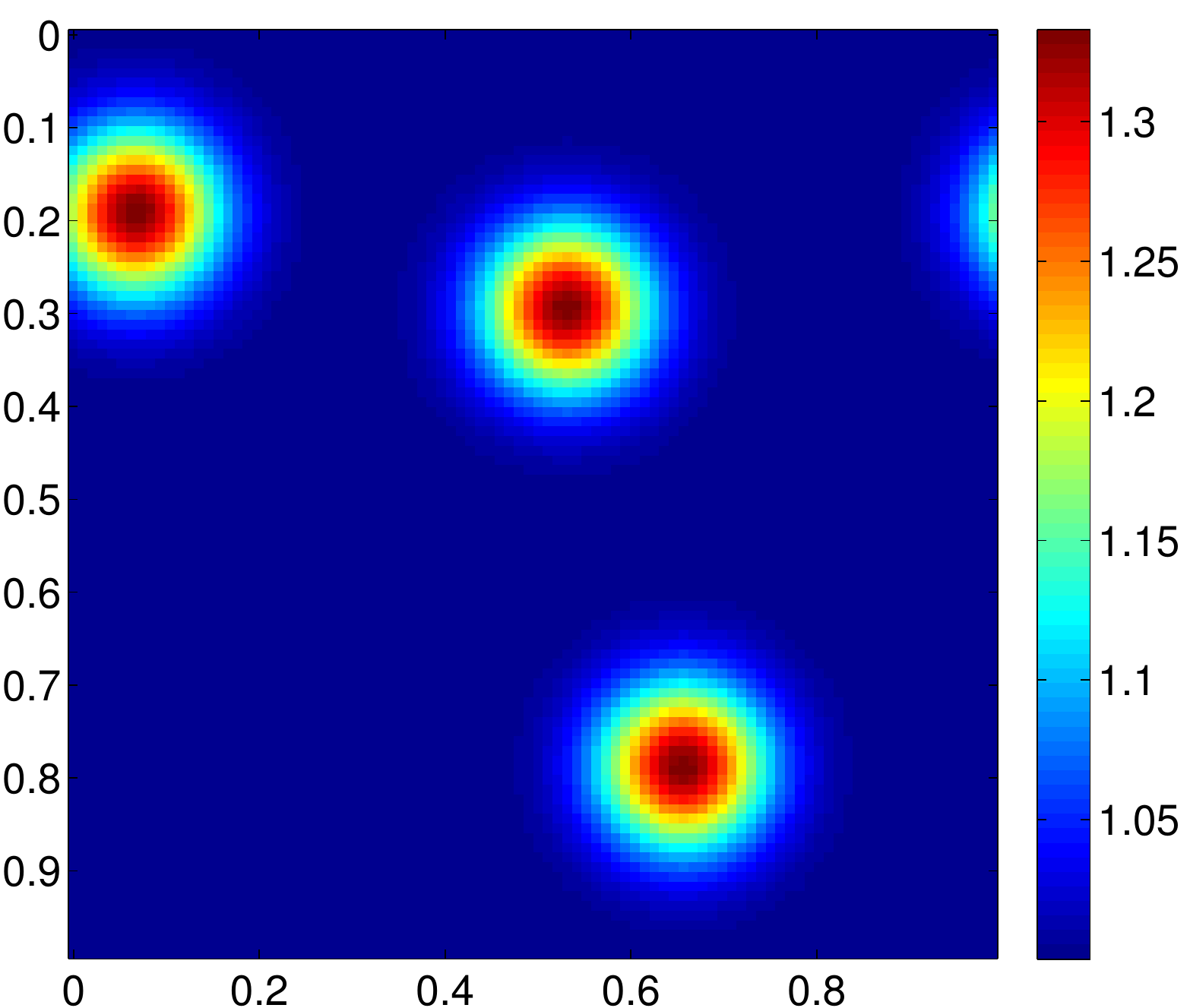} \hspace{0.25in} \includegraphics[height=1.8in]{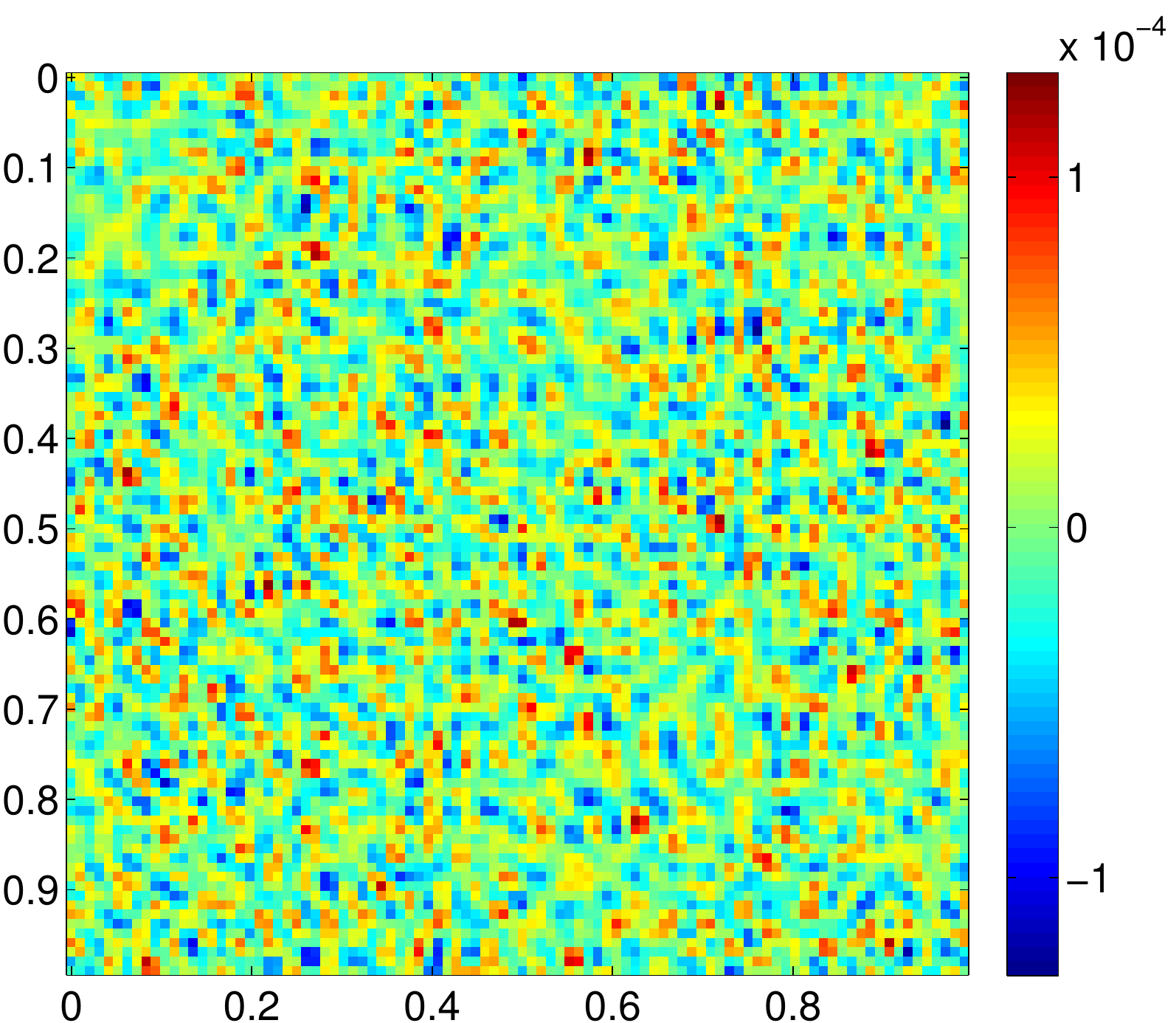}\\
  \end{center}
  \caption{Example 2 of the 3D Helmholtz equation. Top: numerical
    results.  Bottom: $c(x)$ (left) and $u(x)$ (right) for the
    largest $\omega$ value at the middle slice.}
  \label{fig:W32}
\end{figure}

Two similar tests are performed in 3D: (i) in the first one $c(x)$ is
equal to one plus a Gaussian function at the domain, and (ii) in the
second test $c(x)$ is one plus three Gaussians with centers placed
randomly on the middle slice. The right hand side is again a delta
source at the center. The results of these two tests are listed in
Figure \ref{fig:W31} and \ref{fig:W32}.

The numerical results in both 2D and 3D examples show that the
iteration counts remain quite small even for problems at very high
frequency.

\subsection{Schr\"odinger equation}
For the numerical tests of \eqref{eq:schr}, we let $\ell=1/h$ and
consider
\[
(-\Delta + 1/h^2\cdot V(x/h) - 1/h^2\cdot E) u(x) = f(x), \quad x\in \T^d,
\]
where $h=1/n$ is again the step size of the pseudospectral grid. The
energy shift $E$ is chosen to be equal to 2.5 to ensure that there are
about 3 or 4 grid points per wavelength.

\begin{figure}[h!]
  \begin{center}
    \begin{tabular}{|cc|cc|cc|}
      \hline
      $N$ & $b$ & $T_s$(sec) & $T_a$(sec) & $n_p$ & $T_p$(sec) \\
      \hline
      $48^2$ & 3 & 2.1e-01 & 1.3e-02 & 7.0e+00 & 1.0e-01\\
      $96^2$ & 6 & 5.6e-01 & 2.8e-02 & 9.0e+00 & 3.0e-01\\
      $192^2$ & 6 & 2.5e+00 & 7.2e-02 & 2.1e+01 & 2.1e+00\\
      $384^2$ & 12 & 1.3e+01 & 2.1e-01 & 3.2e+01 & 9.6e+00\\
      \hline
    \end{tabular}
    \includegraphics[height=1.8in]{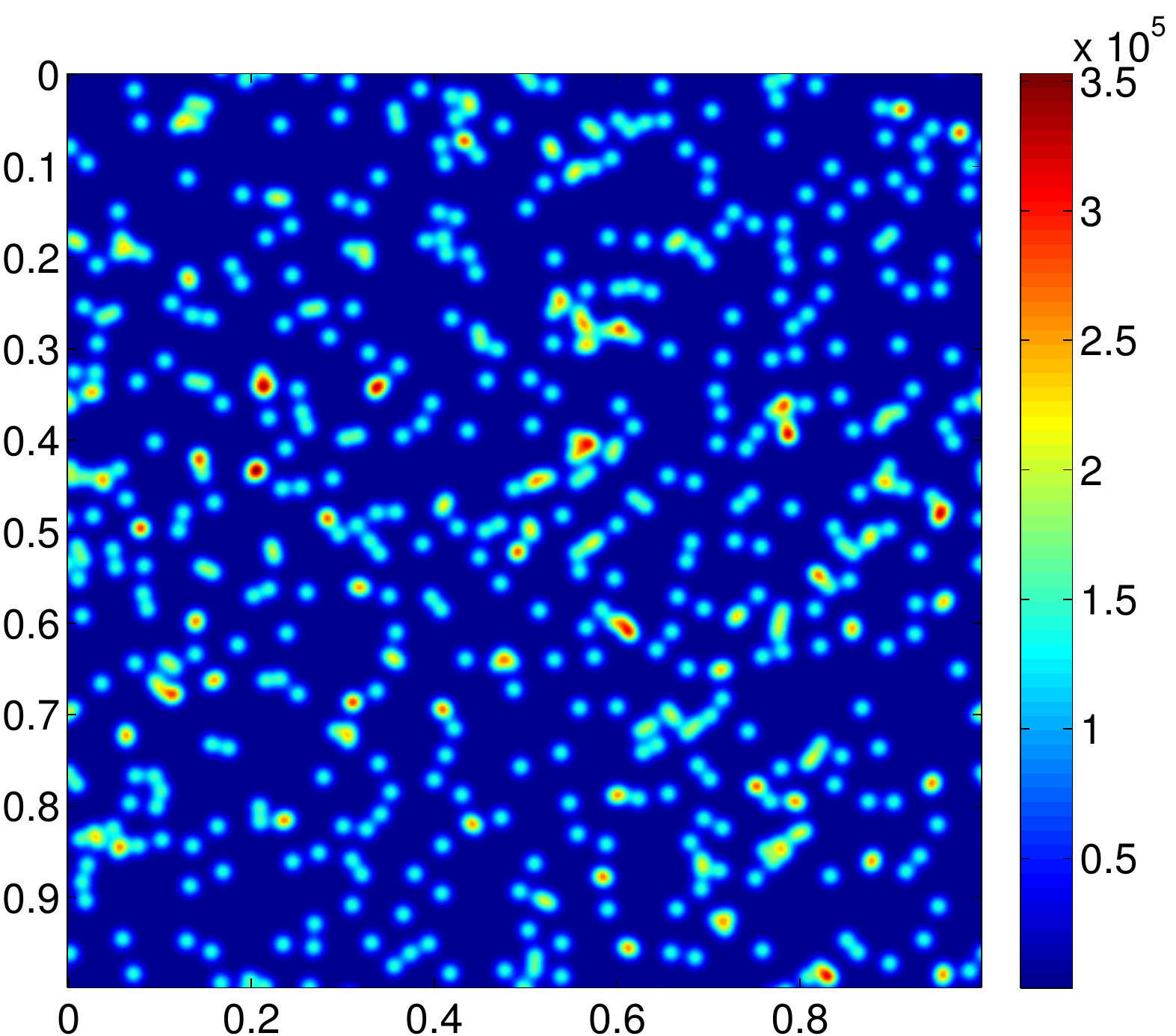} \hspace{0.25in} \includegraphics[height=1.8in]{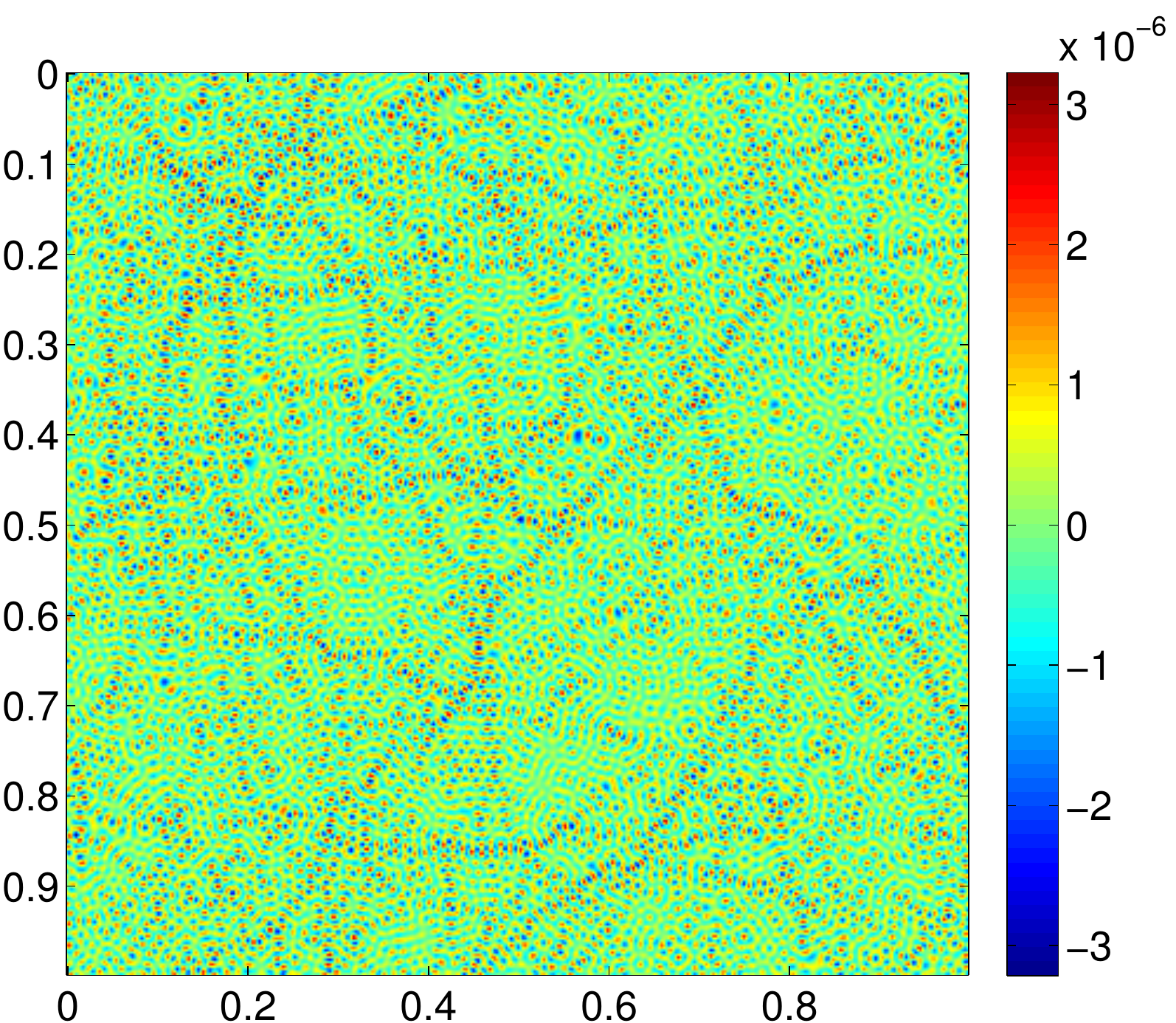}\\
  \end{center}
  \caption{Example 1 of the 2D Schr\"odinger equation. Top: numerical
    results.  Bottom: $1/h^2\cdot V(x/h)$ (left) and $u(x)$ (right) for the
    largest $N$ value.}
  \label{fig:S21}
\end{figure}

\begin{figure}[h!]
  \begin{center}
    \begin{tabular}{|cc|cc|cc|}
      \hline
      $N$ & $b$ & $T_s$(sec) & $T_a$(sec) & $n_p$ & $T_p$(sec) \\
      \hline
      $48^2$ & 3 & 2.2e-01 & 1.3e-02 & 8.0e+00 & 1.1e-01\\
      $96^2$ & 6 & 5.0e-01 & 2.7e-02 & 1.1e+01 & 3.8e-01\\
      $192^2$ & 6 & 3.0e+00 & 7.7e-02 & 2.1e+01 & 2.1e+00\\
      $384^2$ & 12 & 1.2e+01 & 2.9e-01 & 3.8e+01 & 1.3e+01\\
      \hline
    \end{tabular}
    \includegraphics[height=1.8in]{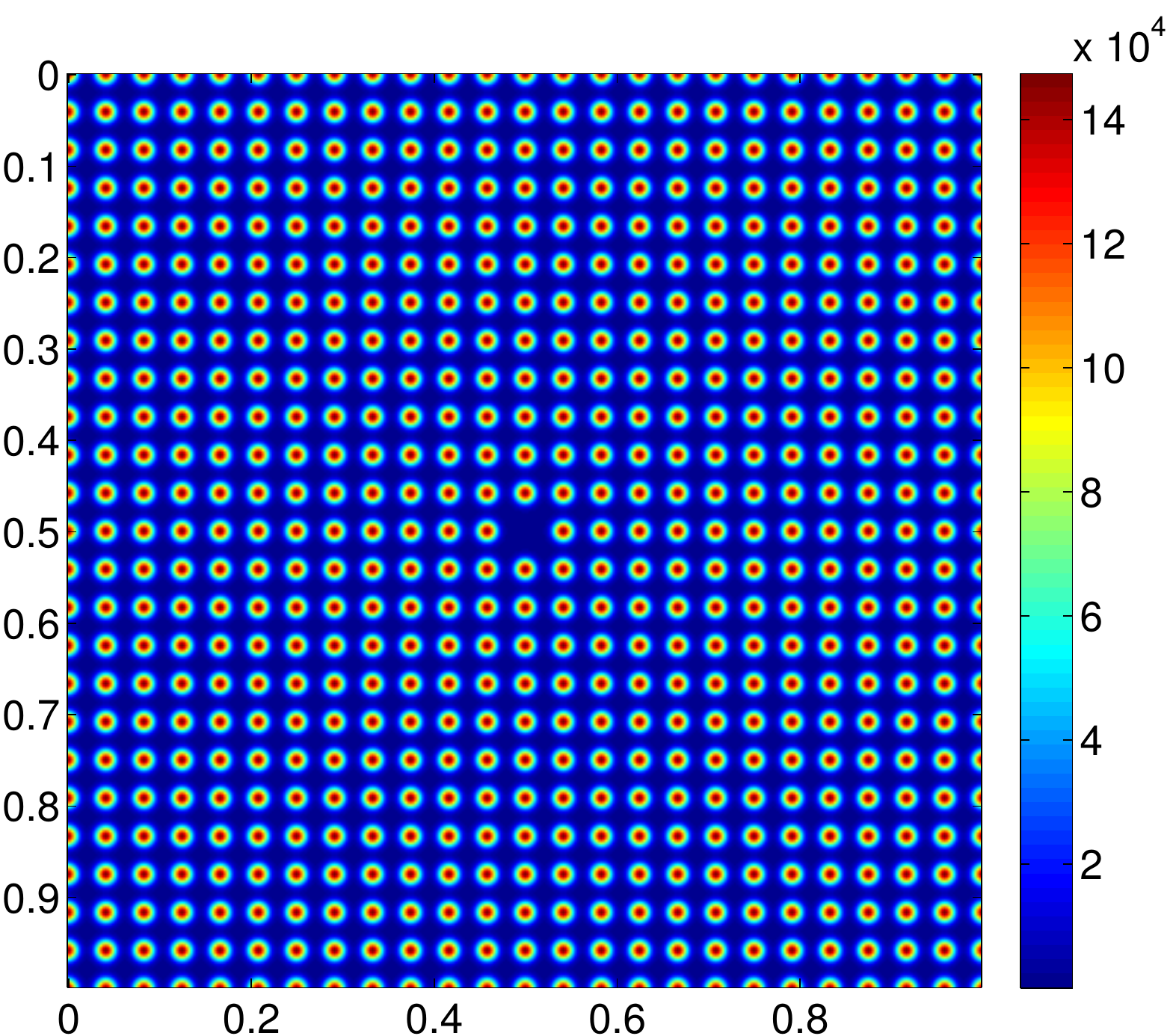} \hspace{0.25in} \includegraphics[height=1.8in]{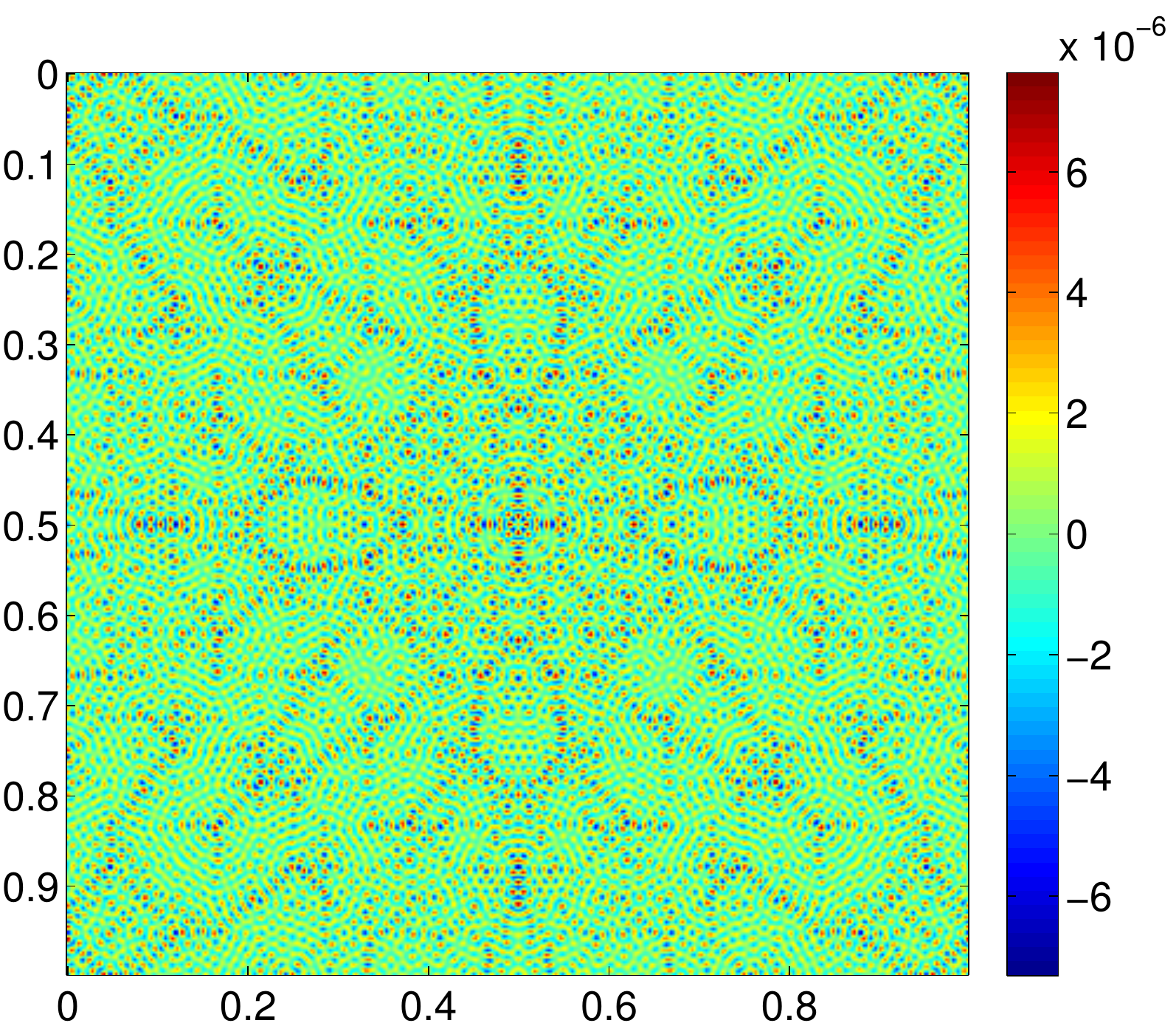}\\
  \end{center}
  \caption{Example 2 of the 2D Schr\"odinger equation. Top: numerical
    results.  Bottom: $1/h^2\cdot V(x/h)$ (left) and $u(x)$ (right) for the
    largest $N$ value.}
  \label{fig:S22}
\end{figure}

Two tests are performed in 2D. In the first one the potential field is
an array of randomly placed Gaussians, while in the second one the
potential field is given by a regular 2D array of Gaussians with one
missing at the center. In both tests, the right hand side is a delta
source at the domain center. The results of these two tests are
summarized in Figures \ref{fig:S21} and \ref{fig:S22}.

\begin{figure}[h!]
  \begin{center}
    \begin{tabular}{|cc|cc|cc|}
      \hline
      $N$ & $b$ & $T_s$(sec) & $T_a$(sec) & $n_p$ & $T_p$(sec) \\
      \hline
      $12^3$ & 3 & 2.9e-01 & 6.6e-03 & 2.0e+00 & 1.1e-02\\
      $24^3$ & 6 & 6.0e+00 & 4.6e-02 & 9.0e+00 & 4.8e-01\\
      $48^3$ & 6 & 1.5e+02 & 4.3e-01 & 1.7e+01 & 8.4e+00\\
      $96^3$ & 12 & 4.5e+03 & 8.4e+00 & 3.1e+01 & 2.8e+02\\
      \hline
    \end{tabular}
    \includegraphics[height=1.8in]{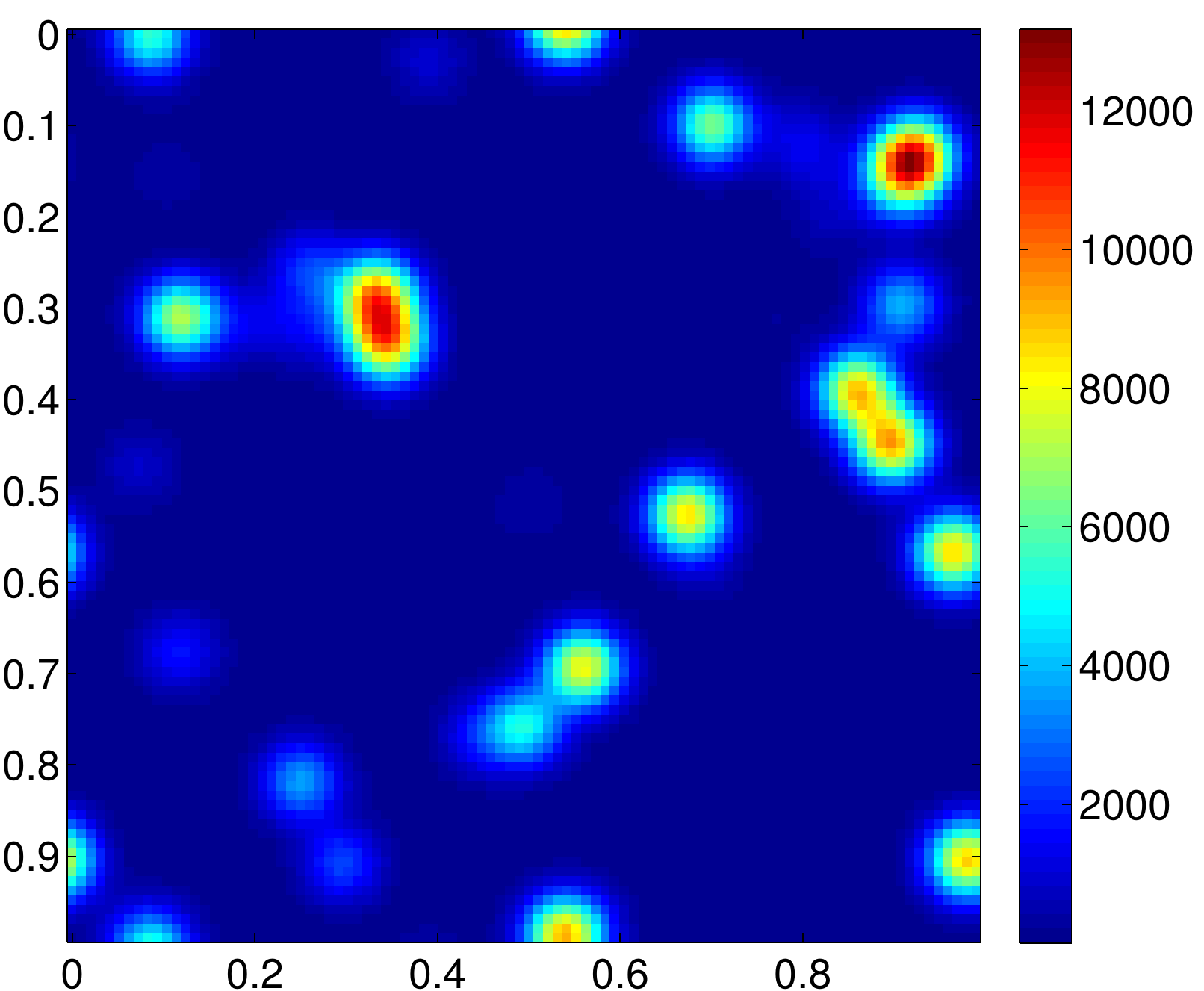} \hspace{0.25in} \includegraphics[height=1.8in]{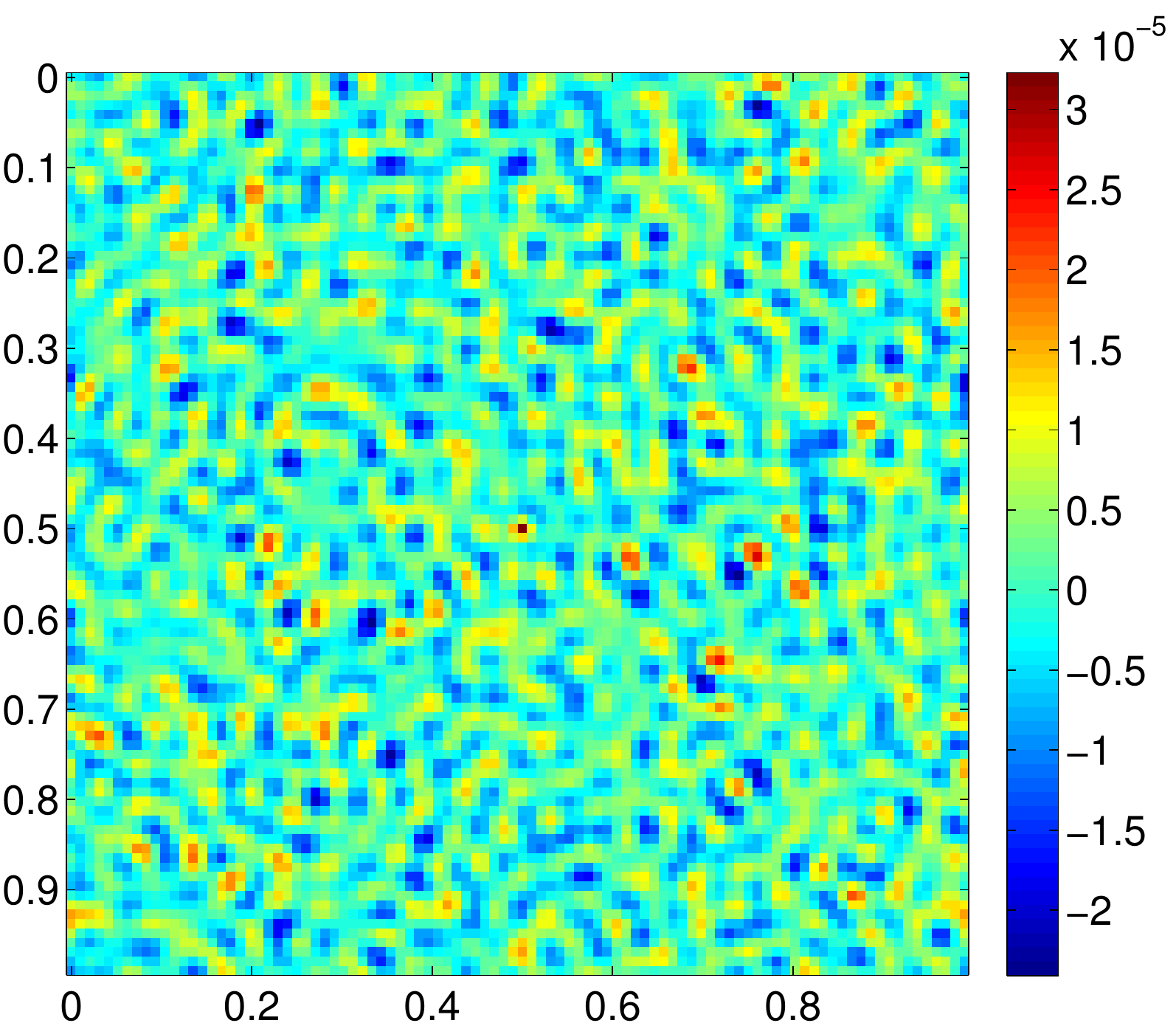}\\
  \end{center}
  \caption{Example 1 of the 3D Schr\"odinger equation. Top: numerical
    results.  Bottom: $1/h^2\cdot V(x/h)$ (left) and $u(x)$ (right) for the
    largest $N$ value at the middle slice.}
  \label{fig:S31}
\end{figure}

\begin{figure}[h!]
  \begin{center}
    \begin{tabular}{|cc|cc|cc|}
      \hline
      $N$ & $b$ & $T_s$(sec) & $T_a$(sec) & $n_p$ & $T_p$(sec) \\
      \hline
      $12^3$ & 3 & 3.4e-01 & 7.1e-03 & 2.0e+00 & 1.6e-02\\
      $24^3$ & 6 & 6.1e+00 & 4.5e-02 & 5.0e+00 & 2.5e-01\\
      $48^3$ & 6 & 1.4e+02 & 5.0e-01 & 6.0e+00 & 3.7e+00\\
      $96^3$ & 12 & 4.6e+03 & 7.4e+00 & 1.2e+01 & 9.9e+01\\
      \hline
    \end{tabular}
    \includegraphics[height=1.8in]{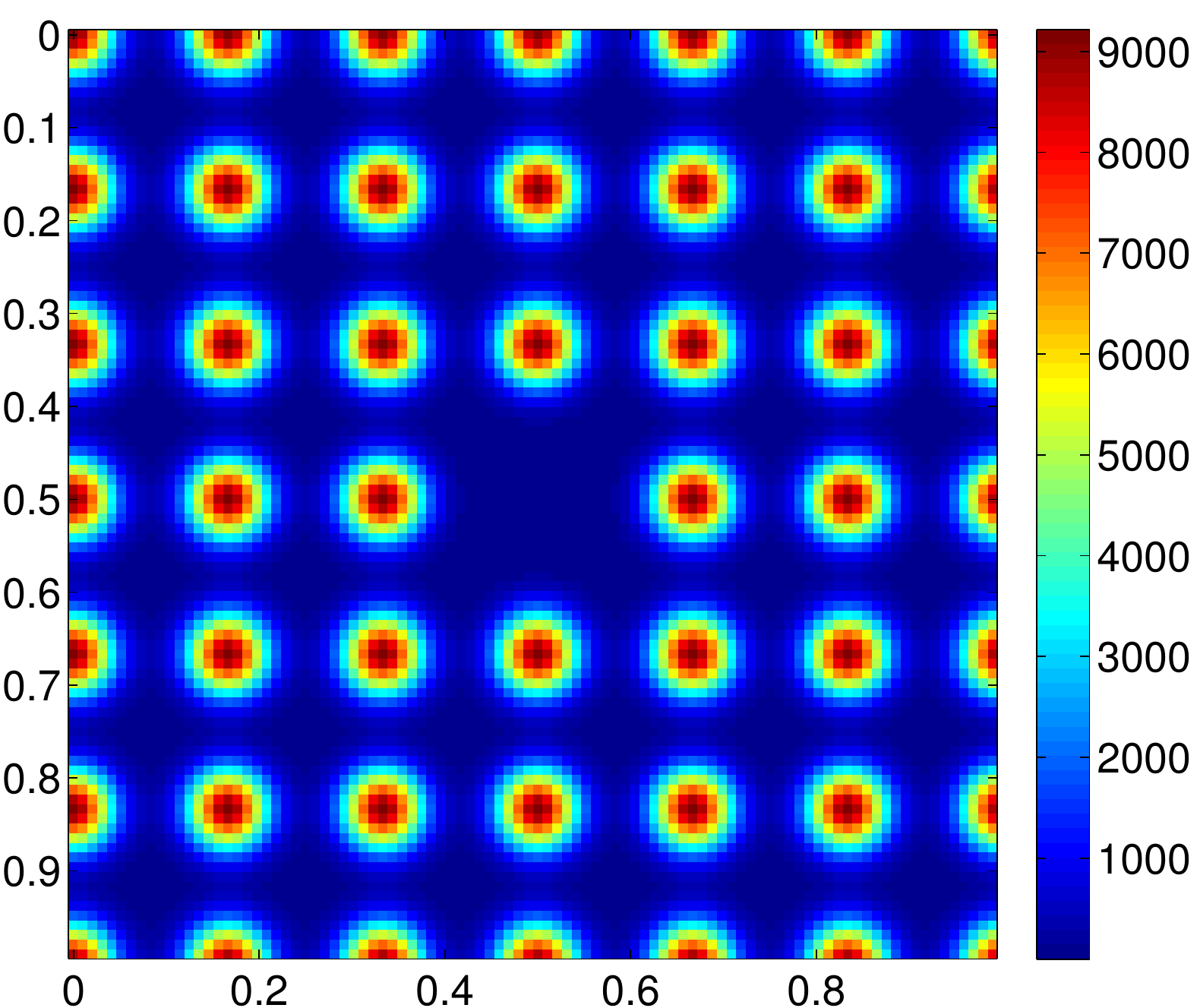} \hspace{0.25in} \includegraphics[height=1.8in]{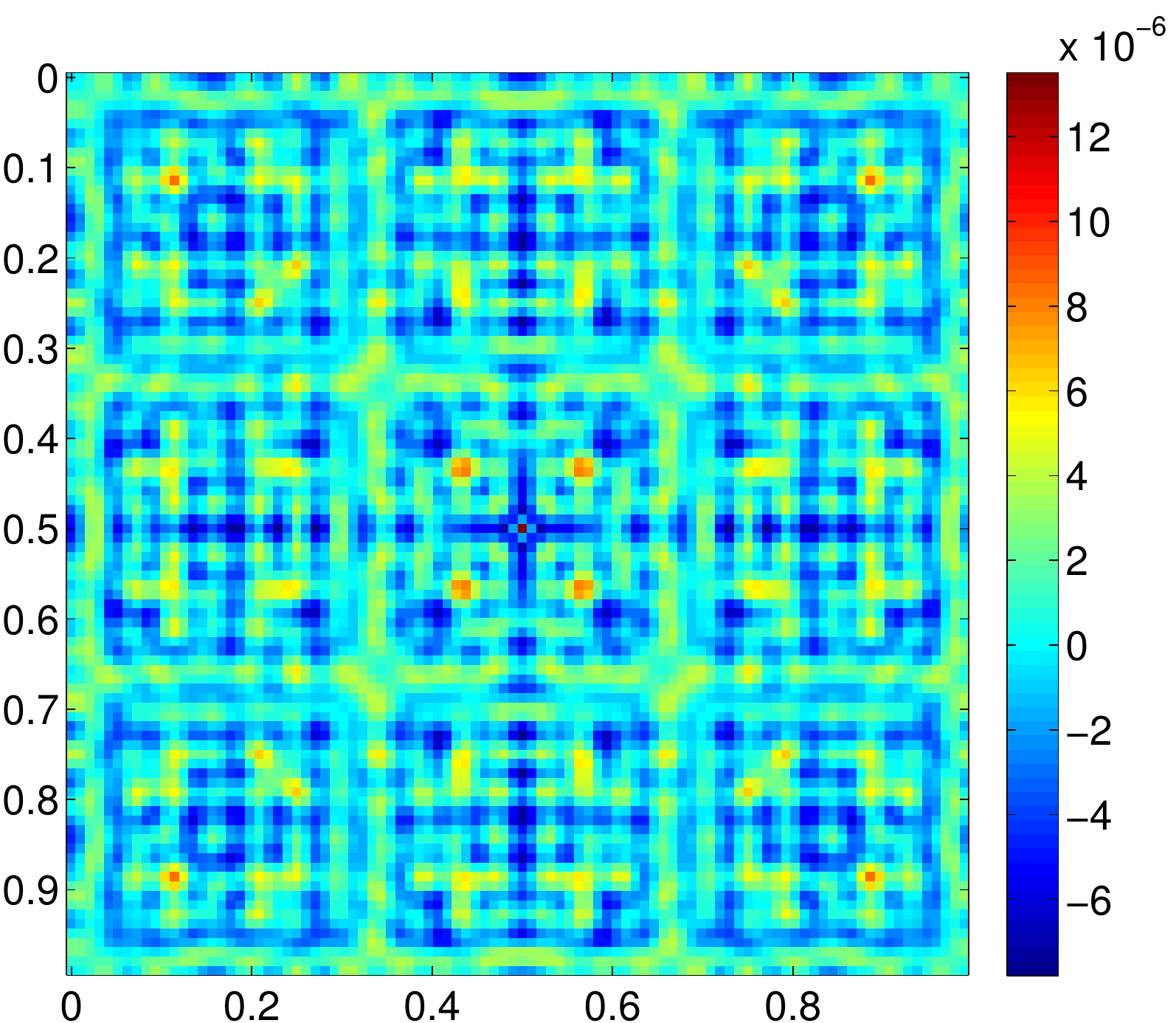}\\
  \end{center}
  \caption{Example 2 of the 3D Schr\"odinger equation. Top: numerical
    results.  Bottom: $1/h^2\cdot V(x/h)$ (left) and $u(x)$ (right) for the
    largest $N$ value at the middle slice.}
  \label{fig:S32}
\end{figure}

Two similar tests are also performed in 3D: (i) in the first one the
potential field is equal to an array of randomly placed Gaussians, and
(ii) in the second test the potential is a regular 3D array of
Gaussians with one missing at the center.  The right hand side is
still a delta source at the domain center. The results of these two
tests are listed in Figure \ref{fig:S31} and \ref{fig:S32}.

In both 2D and 3D examples, the results are qualitatively similar to
the ones of the Helmholtz equation. Though the iteration count
increases with the system size, it remains quite small even for large
scale problems.

%-----------------------------------
\section{Conclusion}

This paper introduces the sparsifying preconditioner for the
pseudospectral approximations of highly indefinite systems on periodic
structures. These systems have important applications in computational
photonics and electronic structure calculation. The main idea of the
preconditioner is to transform the dense system into an integral
equation formulation and introduce a local stencil operator $Q$ for
its sparsification. The resulting approximate equation is then solved
with the nested dissection algorithm and serves as the
preconditioner. This method is easy to implement, efficient, and
results in relatively low iteration counts even for large scale
problems.

In the numerical results, the size $bh$ of the leaf box is of order
$O(n^{1/2} h)$. However, the iteration count still grows roughly
linearly with $\omega$ (for the Helmholtz equation) and with $1/h$
(for the Schr\"odinger equation). An important open question is
whether other sparsity patterns for $Q$ can result in almost frequency
independent iteration counts even for moderate values of $b$.

In computational photonics \cite{joannopoulos-2008}, a more relevant
equation is the Maxwell equation for the electric field $E(x)$:
\[
\left( \nabla \times \nabla \times - \frac{\omega^2}{c^2} \eps(x) \right) E(x) = f(x),
\]
where $\eps(x)$ is the dielectric function. The sparsifying
preconditioner should be extended to this case without much difficulty.

For the density functional theory calculation in computational
chemistry, the Schr\"odinger equation typically has a non-local
pseudopotential term in addition to the local potential term in
\eqref{eq:schr}. An important future work is to extend the sparsifying
preconditioner to address such non-local terms.

The method proposed in this paper provides an efficient way to access
a column or a linear combination of the columns of the Green's
function of the operators in \eqref{eq:helm} and \eqref{eq:schr}. This
can potentially open the door for building efficient and data-sparse
representations of the whole Green's function.

\bibliographystyle{abbrv} \bibliography{ref}

\begin{bibdiv}
\begin{biblist}

\bib{canuto-2007}{book}{
      author={Canuto, C.},
      author={Hussaini, M.~Y.},
      author={Quarteroni, A.},
      author={Zang, T.~A.},
       title={Spectral methods},
      series={Scientific Computation},
   publisher={Springer, Berlin},
        date={2007},
        ISBN={978-3-540-30727-3},
        note={Evolution to complex geometries and applications to fluid
  dynamics},
      review={\MR{2340254 (2009d:76084)}},
}

\bib{duff-1983}{article}{
      author={Duff, I.~S.},
      author={Reid, J.~K.},
       title={The multifrontal solution of indefinite sparse symmetric linear
  equations},
        date={1983},
        ISSN={0098-3500},
     journal={ACM Trans. Math. Software},
      volume={9},
      number={3},
       pages={302\ndash 325},
         url={http://dx.doi.org/10.1145/356044.356047},
      review={\MR{791968 (86k:65030)}},
}

\bib{george-1973}{article}{
      author={George, Alan},
       title={Nested dissection of a regular finite element mesh},
        date={1973},
        ISSN={0036-1429},
     journal={SIAM J. Numer. Anal.},
      volume={10},
       pages={345\ndash 363},
        note={Collection of articles dedicated to the memory of George E.
  Forsythe},
      review={\MR{0388756 (52 \#9590)}},
}

\bib{gottlieb-1977}{book}{
      author={Gottlieb, David},
      author={Orszag, Steven~A.},
       title={Numerical analysis of spectral methods: theory and applications},
   publisher={Society for Industrial and Applied Mathematics, Philadelphia,
  Pa.},
        date={1977},
        note={CBMS-NSF Regional Conference Series in Applied Mathematics, No.
  26},
      review={\MR{0520152 (58 \#24983)}},
}

\bib{joannopoulos-2008}{book}{
      author={Joannopoulos, John~D.},
      author={Johnson, Steven~G.},
      author={Winn, Joshua~N.},
      author={Meade, Robert~D.},
       title={{Photonic Crystals: Molding the Flow of Light (Second Edition)}},
     edition={2},
   publisher={Princeton University Press},
         how={Hardcover},
        date={2008},
        ISBN={0691124566},
  url={http://www.amazon.com/exec/obidos/redirect?tag=citeulike07-20\&path=ASI%
N/0691124566},
}

\bib{orszag-1969}{article}{
      author={Orszag, Steven~A.},
       title={Numerical methods for the simulation of turbulence},
        date={1969},
     journal={Physics of Fluids (1958-1988)},
      volume={12},
      number={12},
       pages={II\ndash 250\ndash II\ndash 257},
  url={http://scitation.aip.org/content/aip/journal/pof1/12/12/10.1063/1.16924%
45},
}

\bib{patera-1984}{article}{
      author={Patera, Anthony~T},
       title={A spectral element method for fluid dynamics: Laminar flow in a
  channel expansion},
        date={1984},
        ISSN={0021-9991},
     journal={Journal of Computational Physics},
      volume={54},
      number={3},
       pages={468 \ndash  488},
  url={http://www.sciencedirect.com/science/article/pii/0021999184901281},
}

\bib{saad-1986}{article}{
      author={Saad, Youcef},
      author={Schultz, Martin~H.},
       title={G{MRES}: a generalized minimal residual algorithm for solving
  nonsymmetric linear systems},
        date={1986},
        ISSN={0196-5204},
     journal={SIAM J. Sci. Statist. Comput.},
      volume={7},
      number={3},
       pages={856\ndash 869},
         url={http://dx.doi.org/10.1137/0907058},
      review={\MR{848568 (87g:65064)}},
}

\bib{trefethen-2000}{book}{
      author={Trefethen, Lloyd~N.},
       title={Spectral methods in {MATLAB}},
      series={Software, Environments, and Tools},
   publisher={Society for Industrial and Applied Mathematics (SIAM),
  Philadelphia, PA},
        date={2000},
      volume={10},
        ISBN={0-89871-465-6},
         url={http://dx.doi.org/10.1137/1.9780898719598},
      review={\MR{1776072 (2001c:65001)}},
}

\bib{ying-2014}{article}{
      author={{Ying}, L.},
       title={{Sparsifying Preconditioner for the Lippmann-Schwinger
  Equation}},
        date={2014-08},
     journal={ArXiv e-prints},
      eprint={1408.4495},
}

\end{biblist}
\end{bibdiv}

\end{document}